\documentclass[10pt,a4paper]{amsart}

\usepackage{graphicx}
\usepackage[a4paper,top=2.5cm,bottom=3.5cm,left=1.7cm,right=2.7cm]{geometry}

    \newcommand{\Q}{{\mathbb Q}}

\newcommand{\hd}{\hat{\delta}}  \newcommand{\hb}{\hat{\beta}}
\newcommand{\of}{\omega_{F}}   \newcommand{\om}{\omega_{m}}

\theoremstyle{definition}
\newtheorem{defn}{Definition}

\theoremstyle{plain}
      
\newtheorem{prop}[defn]{Proposition}
\theoremstyle{remark}

\newtheorem*{nb}{Remark}

\title{Stability Analysis of the Helmholtz Oscillator with Time Varying Mass.}
\author{Richard A. Morrison}
\date{\today}
\begin{document}
\maketitle
\begin{abstract}
We study the effects of periodically time varying mass on the stability of the Helmholtz oscillator, which,
when linearised, takes the form of Ince's equation and exhibits parametric resonance. The resonance regions in the
parameter space are mapped and we use the Melnikov function to demonstrate that the parametric instability does not affect
the qualitative dynamics of the nonlinear system.
\end{abstract}

\section{INTRODUCTION}
The recent activity in the application of Nonlinear Dynamics to problems in Naval
Architecture is exemplified by the papers in \cite{SpyrouThompson2000}. One
particular success has been a more rational understanding of beam
sea capsize and the development of a robust
mathematical toolkit for the study of this problem.

Following the work of Wright and Marshfield, \cite{wm79} it has become common to study the roll dynamics using
a one degree of freedom nonlinear oscillator taking the form
\begin{equation}\ddot{\theta} + g(\dot{\theta}) +G_{Z}(\theta) = F(t).\end{equation}
Here $\theta$ is the roll angle relative to the wave normal,
$g(\dot{\theta})$ is the, often linear--viscous, damping function,
and $G_{Z}(\theta)$ (so denoted as it is usually some polynomial
fit to the $GZ$ curve) encapsulates the restoring moment. The forcing, $F(t)$, is often
taken to be a simple sinusoidal function, see for example \cite{vir89}, although stochastic
models of forcing are also studied, e.g. \cite{jiangetal1996}.

Observing that capsize is equivalent to the roll angle exceeding
the angle of vanishing stability leads to the consideration of the
more general dynamical problem of escape of a particle from a
potential well. The Helmholtz oscillator
\begin{equation}\ddot{x}+\beta\dot{x}+x-x^2=F\sin{\omega t},\label{eq:HT}\end{equation}
known to Naval--Architects as the Helmholtz--Thompson equation, provides a simple archetype to study the phenomena of escape directly.
It has been central to
a number of developments. Thompson's paper \cite{thomp89} provides an overview of the dynamic response of
the system which has been studied experimentally by Gottwald et. al. \cite{gottwaldetal1995}. The inhibition of
chaotic escape is considered in the context of \eqref{eq:HT} in \cite{balibreaetal1998} and a more general approach
to the problem is discussed in \cite{lencirega2000}. The equation finds direct application in the study of
bubble dynamics \cite{kangleal1990} and is much discussed in the Naval Architecture literature, see \cite{thomp97}.
The engineering integrity diagram, \cite{solimanthomspon1989} and the use of Melnikov theory to predict parameter values for which
erosion of basins of attraction takes place were developed in this context. These concepts continue to find fruitful application in quantification
of capsize resistance, see \cite{Spyrouetal2002}.

\bigskip
Providing a reasonable nonlinear dynamical theory of the
stability of a damaged ship remains a considerable challenge. Compared with intact
stability there are a number of dynamical differences which effect
ship motion. The mass of a progressively flooding ship changes
with time. The motion becomes coupled with that of water on deck
or in compartments. The water itself provides non-trivial
modelling problems. Masses of water lead to biases which may also
vary with time.

The increased complexity means that one is forced to make a choice
between sacrificing dynamical information by making simplifying
assumptions or dynamical insight by building less simple comprehensive models.

When making assumptions to simplify a model enough that some
theoretical analysis is possible there is a danger of assuming so
much that all physical meaning is lost. The balance between
finding a model that is simple enough to analyze yet still manages
to capture the essence of some real physical dynamics is difficult
to attain. One successful example is found in Murashige et al's
analysis of flooded ship motion, see \cite{Murashigeetal2000} and
references therein.

In this paper we choose a different approach. The dynamic phenomena listed above will have
varying importance in building a comprehensive theory of damage stability. Rather than start with a
physical model and make simplifying assumptions, we have chosen to put one of these effects, time varying mass,
directly under the microscope. Whilst some properties of time varying mass dynamics which will influence the
behaviour of any model whose mass has a periodic component are discussed, we make no claim to direct physical relevance.

\bigskip

The model we study is a Helmholtz oscillator with time
varying mass. We use the generalised form of Newton's Second Law,
\begin{align} F &=
\frac{d}{dt}\left(m(t)\dot{x}\right)\nonumber\\
&=\dot{m}(t)\dot{x} + m(t) \ddot{x},
\end{align}
so that an oscillator modelled as
$m\ddot{x} + g(x,\dot{x}) = F(t)$ whose mass varies with time becomes
\begin{equation}m(t)\ddot{x}+g(x,\dot{x}) + \dot{m}(t) \dot{x} = F(t).\end{equation}
Equation \eqref{eq:HT} becomes
\begin{equation}
\ddot{x} + \frac{\beta + \dot{m}(t)}{m(t)}\dot{x} + \frac{x-x^2}{m(t)}=F\frac{\sin{\of t}}{m(t)}\label{eq:httvm_general}
\end{equation}
where we have divided through by $m(t)$.

In this paper we study mass varying sinusoidally with time and choose the simple form
\begin{equation} m(t) = \delta + \gamma \sin{\om t}. \label{eq:mass}\end{equation}
Substituting \eqref{eq:mass} into \eqref{eq:httvm_general} and writing $\hd = 1/\delta$ (in order to
simplify the perturbation calculations which follow) we have
\begin{equation}
\ddot{x} + {\hd\frac{\beta+\gamma\om\cos{\om
t}}{1+\gamma\hd\sin{\om t}}\dot{x} +
\hd\frac{x-x^2}{1+\gamma\hd\sin{\om t}}} \\= \hd F\frac{\sin{\of
t}}{1+\gamma\hd\sin{\om t}}.\label{eq:ht_tvm}
\end{equation}

\section{LINEAR STABILITY ANALYSIS OF FIXED POINTS}
Our first step is to consider the linear stability of the two fixed points of the system with $F=0$. We start by scaling
time in order to normalise the frequency of the mass variation. Setting $\tau = \om t$ yields
\begin{equation}
\ddot{x} + \hd\frac{\frac{\beta}{\om}+\gamma\cos{t}}{1+\gamma\hd\sin{t}}\dot{x}+\frac{\hd}{\om^2}\frac{x-x^2}{1+\gamma\hd\sin{t}}=0.\label{eq:httmvnd}
\end{equation}
after relabelling $t = \tau$. Writing \eqref{eq:httmvnd} as a system
\begin{subequations}
\begin{align}
\dot{x} &= y \\
\dot{y} &= -\frac{\hd}{\om^2}\frac{x-x^2}{1+\gamma\hd\sin{t}}-\hd\frac{\frac{\beta}{\om}+\gamma\cos{t}}{1+\gamma\hd\sin{t}}y
\end{align}
\end{subequations}
we see that the two fixed points are $(x,y) = (0,0)$ and $(1,0)$. We know that for $\gamma = 0$ the origin is a stable focus or
centre as $\beta>0$ or $\beta=0$ respectively. The point $(1,0)$ is a saddle independent of the value of $\beta$.

For both fixed points the linearised system takes the form
$$ \dot{\textbf{x}} = A(t) \textbf{x}$$ where $A(t)=A(t+2\pi)$ is a periodic matrix and we are in the setting of
Floquet theory \cite{js7799}.
For the saddle point we have
\begin{equation}
A_{(1,0)}(t) =\left(\begin{array}{cc}0&1\\ \frac{\hd}{1+\gamma\hd\sin{t}} &-\hd\frac{\frac{\beta}{\om}+\gamma\cos{t}}{1+\gamma\hd\sin{t}} \end{array}\right).\end{equation}
The origin gives
\begin{equation}
A_{(0,0)}(t) =\left(\begin{array}{cc}0&1\\ -\frac{\hd}{1+\gamma\hd\sin{t}} &-\hd\frac{\frac{\beta}{\om}+\gamma\cos{t}}{1+\gamma\hd\sin{t}} \end{array}\right).\label{eq:dforigin}\end{equation}

\subsection*{THE FIXED POINT AT $(1,0)$}
The linear stability analysis of this fixed point does not require the full machinery of Floquet theory.
The ODE linearised about this point is
$$\ddot{x}+\frac{\frac{\beta}{\om}+\hd\gamma\cos{t}}{1+\gamma\hd\sin{t}}\dot{x}-\frac{\hd}{\om^2\left(1+\gamma\hd\sin{t}\right)}x=0.$$
Upon expanding $x=x_0+\gamma x_1 + O(\gamma^2)$ and $\hd=\hd_0+\gamma\hd_1+O(\gamma^2)$ and collecting terms in $\gamma$
we see that the $O(1)$ equation is
$$\ddot{x_0}+\frac{\beta}{\om}\dot{x}_0-\hd_0 x_0=0.$$ The origin in this equation is stable only when $\hd_0<0$ which is
physically meaningless. We conclude that the fixed point of the
nonlinear system is unstable for all $\hd > 0$, regardless of the
value of  $\beta$ and $\om$.

\subsection*{THE FIXED POINT AT THE ORIGIN}
The ODE linearised about the origin is
\begin{equation}\ddot{x}+\frac{\frac{\beta}{\om}+\hd\gamma\cos{t}}{1+\gamma\hd\sin{t}}\dot{x}+\frac{\hd}{\om^2\left(1+\gamma\hd\sin{t}\right)}x=0\label{eq:linearorigin}.\end{equation}
and in this case there are a number of regions where a Mathieu type instability
occurs. It follows from Floquet theory that transition curves in the parameter space correspond to
solutions of the system which have periods equal to and double the fundamental period of the system. Thinking about
periodicity is easier when the fundamental periods under consideration are $2\pi$ and $4\pi$ so we scale time again setting
$t=2\tau.$ As we wish to follow Rand's treatment of Ince's equation
\begin{equation*}
(1+a\cos{2t})\ddot{x}+b\sin{2t}\dot{x}+ (c+d\cos{2t})x=0,
\end{equation*}
see \cite{rand2002}, we also translate time $\tau\mapsto\tau+2\pi$. Relabelling $\tau=t$ for
notational convenience we have
\begin{equation}
(1+\hd\gamma\cos{2t})\ddot{x}+2\left(\frac{\beta}{\om}-\hd\gamma\sin{2t}\right)\dot{x}+
4\frac{\hd}{\om^2}x=0.\label{eq:linearorigin_ince}
\end{equation}
When $\beta=0$ this is Ince's equation with $a=\hd\gamma$, $b=-2\hd\gamma$, $c=4\frac{\hd}{\om^2}$ and $d=0$.
We will consider $\beta\neq 0$ when the $\beta=0$ dynamics are understood.
\begin{nb}
The scaling of time has no effect on the stability of the origin in the linearised system. The stable parameters for
system (\ref{eq:linearorigin}) are exactly those that are stable for system (\ref{eq:linearorigin_ince}).
\end{nb}

As transition curves are characterised by the existence of solutions of period $2\pi$ and $4\pi$ we use harmonic balance to
seek solutions of this form. Substituting
$
x=\sum_{n=0}^{\infty}a_n\cos{nt}+b_n\sin{nt}
$
into \eqref{eq:linearorigin_ince}, combining trigonometric terms and equating coefficients of $\cos{nt}$ and $\sin{nt}$
leads to an infinite set of linear homogeneous equations in the Fourier coefficients $a_i$ and $b_i$. The first few of these
are

\begin{gather*}
4\frac{\hd}{\om^2}a_0 =0\\
\left(1+\frac{\hd\gamma}{2}-4\frac{\hd}{\om^2}\right)a_1-\frac{3}{2}\hd\gamma a_3=0\\
\left(1-\frac{\hd\gamma}{2}-4\frac{\hd}{\om^2}\right)b_1-\frac{3}{2}\hd\gamma b_3=0\\
\left(4\frac{\hd}{\om^2}-4\right)a_2 + -4\hd\gamma a_4 =0\\
\vdots\nonumber
\end{gather*}

The existence of a transition curve in the parameter space is equivalent to the existence of a nontrivial solution to this system of equations.
This in turn is equivalent to the vanishing of the associated determinants of the system. This last condition says that the
transition curves in the parameter space are composed of points which satisfy the equations
\begin{subequations}
\begin{align}
\triangle_{a_{2i}} &=
\left|\begin{array}{cccccc}
4\frac{\hd}{\om^2}&-2\hd\gamma&0&0&0&\cdots\\
0 & 4\frac{\hd}{\om^2}-4&-4\gamma\hd&0&0&\cdots\\
0 & -4\gamma\hd & 4\frac{\hd}{\om^2}-16&-12\gamma\hd&0&\cdots\\
0 & 0 & -12\gamma\hd & 4\frac{\hd}{\om^2}-36 & -24\gamma\hd & \cdots\\
\vdots&\vdots&\vdots&\vdots&\vdots&\ddots
\end{array}\right| \ =0\label{eq:detaeven=0},\\
\triangle_{b_{2i}}&=
\left|
\begin{array}{cccccc}
4\frac{\hd}{\om^2}-4&-4\gamma\hd&0&0&0&\cdots\\
-4\gamma\hd & 4\frac{\hd}{\om^2}-16&-12\gamma\hd&0&0&\cdots\\
0 & -12\gamma\hd & 4\frac{\hd}{\om^2}-36 & -24\gamma\hd &0& \cdots\\
0&0& -24\gamma\hd & 4\frac{\hd}{\om^2}-64& -40\gamma\hd & \cdots\\
\vdots&\vdots&\vdots&\vdots&\vdots&\ddots
\end{array}
\right| \ =0,\label{eq:detbeven=0}\\
\triangle_{a_{2i+1}}&=
\left|
\begin{array}{cccccc}
4\frac{\hd}{\om^2}-1+\frac{\gamma\hd}{2}&-\frac{3}{2}\gamma\hd&0&0&0&\cdots\\
-\frac{3}{2}\gamma\hd&4\frac{\hd}{\om^2}-9&-\frac{15}{2}\gamma\hd&0&0&\cdots\\
0&-\frac{15}{2}\gamma\hd&4\frac{\hd}{\om^2}-25&-\frac{35}{2}\gamma\hd&0&\cdots\\
0&0&-\frac{35}{2}\gamma\hd&4\frac{\hd}{\om^2}-49&-\frac{63}{2}\gamma\hd&\cdots\\
\vdots&\vdots&\vdots&\vdots&\vdots&\ddots
\end{array}
\right| \ =0\label{eq:detaodd=0}\\
\intertext{and}
\triangle_{b_{2i+1}}&=
\left|
\begin{array}{cccccc}
4\frac{\hd}{\om^2}-1-\frac{\gamma\hd}{2}&-\frac{3}{2}\gamma\hd&0&0&0&\cdots\\
-\frac{3}{2}\gamma\hd&4\frac{\hd}{\om^2}-9&-\frac{15}{2}\gamma\hd&0&0&\cdots\\
0&-\frac{15}{2}\gamma\hd&4\frac{\hd}{\om^2}-25&-\frac{35}{2}\gamma\hd&0&\cdots\\
0&0&-\frac{35}{2}\gamma\hd&4\frac{\hd}{\om^2}-49&-\frac{63}{2}\gamma\hd&\cdots\\
\vdots&\vdots&\vdots&\vdots&\vdots&\ddots
\end{array}
\right| \ =0.\label{eq:detbodd=0}
\end{align}
\end{subequations}

We note that
$$\triangle_{a_{2i}}=4\frac{\hd}{\om^2}\triangle_{b_{2i}}.$$ This fact has two consequences. First of all it is
clear that $\hd=0$ is a solution of (\ref{eq:detaeven=0}) so that the first transition curve we find is
\begin{equation}
\hd=0.
\end{equation}
Secondly all other zeros of $\triangle_{a_{2i}}=0$ will coincide with the zeros of $\triangle_{b_{2i}}=0$. Upon setting
$\gamma=0$ it is clear that these curves emanate from the points
$$\hd=\left(\frac{2n\om}{2}\right)^2, \ n \in \{1,2,3,\dots\}.$$ However because they coincide the regions they bound have zero width and no instability occurs. The
coincidence of two transition curves and associated vanishing of a regime of instability is known as \emph{coexistence} see
\cite{magnuswinkler66} or \cite{rand2002} for details.
\begin{nb} The lack of instability is not a consequence of the vanishing width of the unstable regime. Periodic
solutions exist on the transition curves and those on the $\triangle_{a_{2i}}=0$ are even functions of $t$. On the
$\triangle_{b_{2i}}=0$ the solutions are odd. As odd and even solutions are linearly independent we have, on
two coinciding curves, a spanning set for the vector space of solutions. A solution space spanned by bounded
solutions cannot contain an unbounded solution.
\end{nb}

In order to evaluate the solution curves of $\triangle_{a_{2i+1}}=0$ and $\triangle_{b_{2i+1}}=0$ we substitute
$$\hd=\hd_0 + \hd_1\gamma+\hd_2\gamma^2+\hd_3\gamma^3+\hd_4\gamma^4+ \dots$$ into a finite, $n\times n$ where
$n$ is a small integer, truncation of the determinants $\triangle_{a_{2i+1}}$ and $\triangle_{b_{2i+1}}$,
equate to zero and collect coefficients of $\gamma$ in the usual fashion. Using Maple and determinants
of order $25\times 25$ we can easily write down the first few curves.
\begin{figure}
\begin{center}
\includegraphics[width=0.6\linewidth]{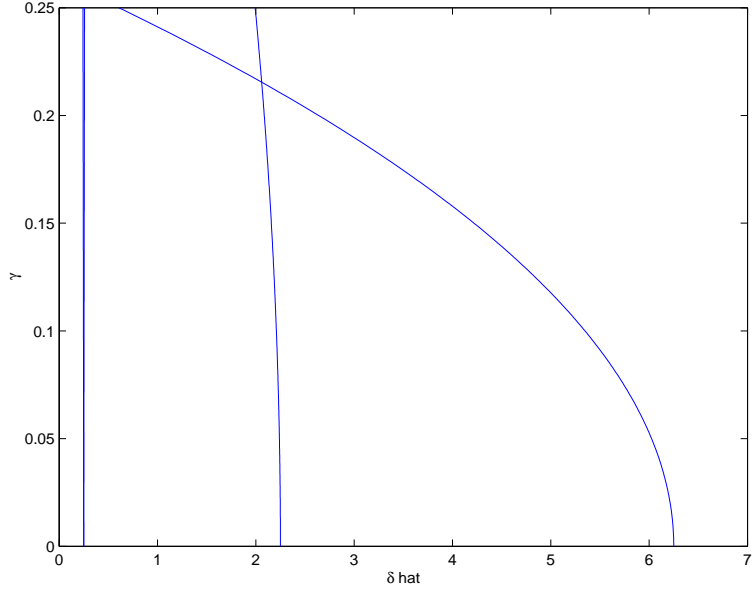}
\end{center}
\caption{The first three regions of parametric instability.}\label{fig:three_regions}
\end{figure}
\begin{subequations}\label{eq:transition_approx}
\begin{align}
\hd&=\frac{\om^{2}}{4}-\frac{\om^{4}}{32}\gamma+O(\gamma^2)\\
\hd&=\frac{\om^{2}}{4}+\frac{\om^{4}}{32}\gamma+O(\gamma^2)\\
\hd&=\frac{9}{4}\om^{2}-{\frac {16767}{4096}}\,\om^{6}{\gamma}^{2}-{
\frac {6561}{131072}}\,\om^{8}{\gamma}^{3}+O({\gamma}^{4})\\
\hd&=\frac{9}{4}\om^{2}-{\frac {16767}{4096}}\,\om^{6}{\gamma}^{2}+{
\frac {6561}{131072}}\,\om^{8}{\gamma}^{3}+O({\gamma}^{4})\\
\hd&={\frac {25}{4}}\,\om^{2}-{\frac {1109375}{12288}}\,{\omega
}^{6}{\gamma}^{2}+{\frac {3030048828125}{1811939328}}\,\om^{10}{
\gamma}^{4}+{\frac {2197265625}{2147483648}}\,\om^{12}{\gamma}^{5
}+O({\gamma}^{6})\\
\hd&={\frac {25}{4}}\,\om^{2}-{\frac {1109375}{12288}}\,{\omega
}^{6}{\gamma}^{2}+{\frac {3030048828125}{1811939328}}\,\om^{10}{
\gamma}^{4}-{\frac {2197265625}{2147483648}}\,\om^{12}{\gamma}^{5
}+O({\gamma}^{6})
\end{align}
\end{subequations}
As the equations of the curves suggest, it is possible to prove:

\begin{prop}\label{prop:transition_split}
The pair of transition curves emanating from the point
$$\hd=\left(\frac{(2k-1)\om}{2}\right)^2$$ coincide to order $\gamma^{2k-2}$. They split at order $\gamma^{2k-1}$.
\end{prop}

\begin{figure}
\begin{center}
\includegraphics[width=0.5\linewidth]{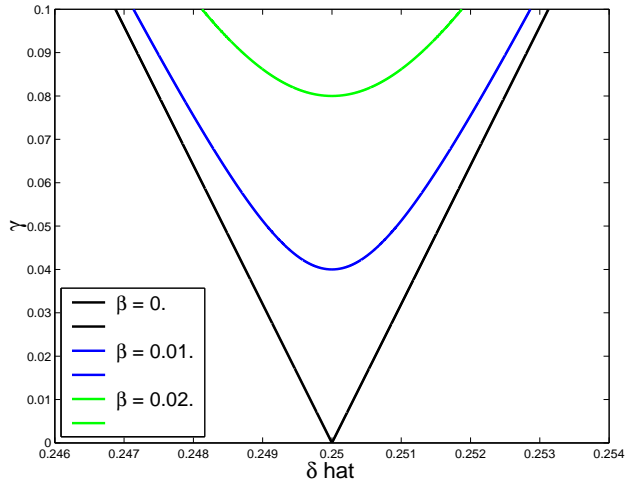}
\end{center}
\caption{The First Region of Instability (With and Without Damping). $\om=1$}\label{fig:beta}
\end{figure}

The geometry of this fact is that as $\hd$ grows the curves becomes thinner and thinner. The areas of the parameter space
in which instability occurs decreases.  Figure \ref{fig:three_regions} illustrates the first three regions of instability,
all of which appear thin on this scale. A closer look at the first region (along with the effect of $\beta>0$, discussed below)
is seen in Figure \ref{fig:beta}.

\begin{figure}
\begin{center}
\includegraphics[width=0.5\linewidth]{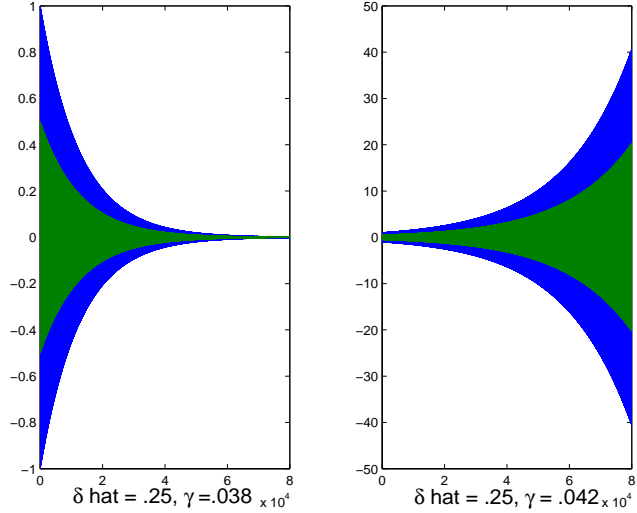}
\end{center}
\caption{Trajectories of the Linear System. The left hand solution has parameters just outside the resonant region. The
right hand, just inside. In both cases the solution were plotted using the Matlab ode45 routine with initial condition $(1,0)$.}\label{fig:beta_traj}
\end{figure}

\begin{nb}
We have not mentioned that a mathematically rigourous treatment would have to address the divergence
of the determinants in their current form. The solution
to this problem involves noting that the individual coefficients of $\gamma^n$ are defined as ratios of infinite
determinants. One must then consider the coefficient as the ratio of two $k\times k$ determinants, let
$k\rightarrow\infty$, and show that the expression converges. This, along with the proof of Proposition
\ref{prop:transition_split}, are essentially exercises in algebra and we do not reproduce the calculations
here.
\end{nb}
When $\beta\neq0$ the above analysis fails as the systems of equations for the coefficients of the Fourier series do not
decouple in a nice way. We use a perturbation argument instead to illustrate the effect of $\beta$. It is necessary
to assume that $\beta$ is of the same order as $\gamma$ which we consider to be a small parameter. Thus we write
$\beta=\gamma\hb$ and we are studying
\begin{equation}\ddot{x}+\gamma\hd\frac{\frac{\hb}{\om}+\cos{t}}{1+\gamma\hd\sin{t}}\dot{x}+\frac{\hd}{\om^2\left(1+\gamma\hd\sin{t}\right)}x=0.\label{eq:perturbation1}\end{equation}
We expand $\hd=\hd_0+\hd_1\gamma+\hd_2\gamma^2+O(\gamma^3)$, $x=x_0+x_1\gamma+x_2\gamma^2+O(\gamma^3)$, substitute and
collect coefficients in the usual way. Considering $\hd_0=\om^2/4$ we see that $$\hd_1 = \pm \frac{\om^3}{32}\sqrt{\om^2-16\hb^2}.$$
Now remembering that $\gamma\hb=\beta$ we have
\begin{align}
\hd&=\frac{\om^2}{4} \pm \frac{\om^3}{32}\sqrt{\om^2-16\hb^2}\gamma + O(\gamma^2)\nonumber\\
&=\frac{\om^2}{4}\pm \frac{\om^3}{32}\sqrt{\gamma^2\om^2-16\beta^2} + O(\gamma^2).
\end{align}
In Figure \ref{fig:beta} we see that the effect of $\beta>0$ is to lift the curves from the $\gamma=0$ axis,
increasing the amplitude of mass variation that is required to produce an instability. Figure \ref{fig:beta_traj} illustrates the numerical response of the system in $(x,t)$ space. The left hand graph
corresponds to parameters outwith the region of parametric resonance, the right has parameters inside the region.

\section{THE MELNIKOV FUNCTION FOR THE NONLINEAR SYSTEM}
Having understood the linear dynamics of the two fixed points we now wish to consider the full nonlinear problem. We will
use Melnikov Theory (see \cite{thomp96} or \cite{wig90}) to understand which parameter values support homoclinic tangling
(and therefor safe basin erosion) and compare the analytic results with numerically calculated integrity curves.

To calculate the Melnikov Function we need to consider \eqref{eq:ht_tvm} as a nearly integrable system taking the form
\begin{equation}\dot{\textbf{x}} = \textbf{f}(\textbf{x}) + \epsilon \textbf{g}(\textbf{x},t).\label{eq:mref}\end{equation}
We assume that $F,\gamma,$ and $\beta$
are $O(\epsilon)$ and expand $$\frac{\hd}{1+\gamma\hd\sin{\om t}} = 1-\gamma\hd\sin{\om t} + O(\gamma^2).$$ Writing
$F = \epsilon F$, $\gamma=\epsilon\gamma$ and $\beta=\epsilon\beta$ for notational convenience we have \eqref{eq:ht_tvm}
in the form \eqref{eq:mref} and the Melnikov function is
\begin{multline}
M(t_0) = \int_{-\infty}^{\infty} \bigg(F\hd\sin{\left(\of (t+t_0)\right)} \ y_h(t)
-{\hd\beta}y_h^2(t)- \hd\gamma\om\cos{\left(\om(t+t_0)\right)} \ y_h^2(t)\\+{\gamma\hd^2}\sin{\left(\om(t+t_0)\right)}(x_h(t)-x_h^2(t))y_h(t)\bigg)dt.\nonumber
\end{multline}
Here
$$x_h(t) = 1-\frac{3}{1+\cosh{\sqrt{\hd}t}}, \> y_h = \dot{x}_h$$ is the parameterised orbit homoclinic to $(1,0)$.
Using trigonometric reduction, changes of variables and residue calculus this evaluates as
\begin{equation}
M(t_0) = \frac{6F\of^2\pi\cos{\of t_0}}{\sinh{\frac{\pi\of}{\sqrt{\hd}}}} - \frac{6}{5}\hd^{\frac{3}{2}}\beta  - \gamma\frac{3}{5}\frac{\pi\om^2}{\hd}\frac{\cos{\om t_0}}{\sinh{\frac{\pi\om}{\sqrt{\hd}}}}\left(\hd^2-\om^4\right).\label{eq:Melnikov}\end{equation}
Notice that when $\hd=1$ and $\gamma=0$ this is identical to the Melnikov function for the standard Helmholtz oscillator, \cite{thomp96}.
\begin{nb} In cases where $\om\neq\of$ and in particular when $\om/\of \not\in\Q$ we must rely on the
generalisation of Melnikov theory to quasi--periodic perturbations due to Meyer and Sell \cite{meyersell89} and Wiggins \cite{wig88}.
In this paper \emph{Melnikov Function} is used as a synonym for \emph{Generalised Quasi--Periodic Melnikov Function}.
\end{nb}

\subsubsection*{The Effect of Linear Resonance on Basin Erosion}
\begin{figure}
\begin{center}
\includegraphics[width=0.6\linewidth]{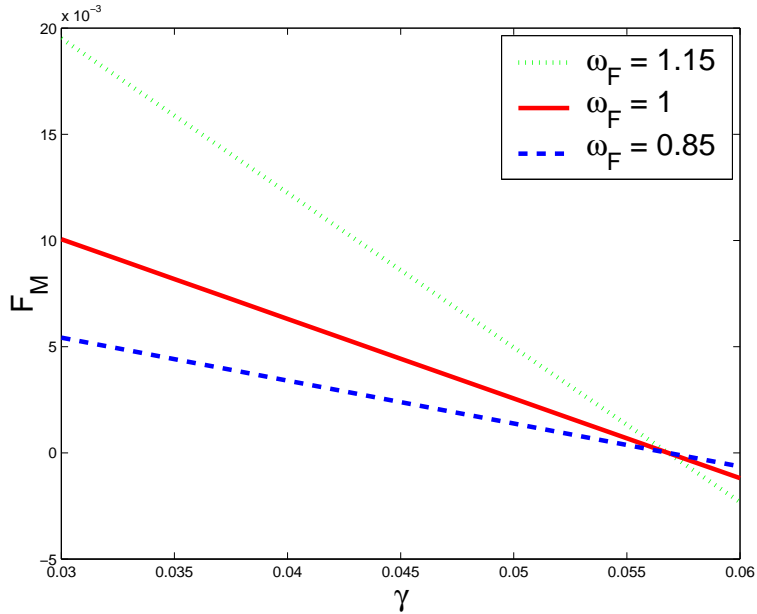}
\end{center}
\caption{$F_M$ as a function of  $\gamma$. Even as $\gamma$ passes over the stability boundary at
$\gamma = 0.04$ the function remains linear.}\label{fig:Fgamma}
\end{figure}

To estimate the minimal $F$ for which erosion is expected we set $M(t_0)=0$ and solve
the function for $F$, selecting $t_0$ to minimise as appropriate. The first question we
wish to answer is does basin erosion differ depending on whither or not we are located within a
region of parametric instability. Consider the region illustrated in Figure \ref{fig:beta}. In this
case $\om=1$ and the region is centred at $\hd=0.25$. First we fix $\beta=0.01$, $\hd = 0.25$ and vary $\gamma$ to
cross the boundary of the region. With these parameters the Melnikov estimate of $F$ looks like
$$F_M(\gamma,\of) = \frac{\left(.0015\sinh{2\pi}-2.25\pi\gamma\right)\sinh{2\pi\of}}{6\pi\of^2\sinh{2\pi}}$$
and for fixed $\of$ this is just a linear function of $\gamma$, see Figure \ref{fig:Fgamma}.

If the presence of
parametric instability had any interesting effect upon basin erosion we would expect some effect on $F_M$
as $\gamma$ moved across the stability boundary. There is, however, no cusp, blip,
change in gradient or other effect in the Melnikov estimate as the parameters move across the transition curve.

The parametric resonance has no effect on the speed of
homoclinic tangling and therefor no import for basin erosion. As the tangling and erosion processes are governed by the
invariant manifolds of the saddle this result is unsurprising. Intuitively the homoclinic orbit \emph{does not see} the effect of
the local dynamics at the other fixed point. Time varying mass is, however, important in its quantitative effects on
basin erosion.
\subsubsection*{Time Varying Mass as Excitation}
Notice in Figure \ref{fig:Fgamma} that when $\gamma$ is slightly greater than $0.055$ the lines all cross $F_M=0$. To understand this we consider the Melnikov
Function with $F=0$ as an estimator of $\gamma$. We have
\begin{equation*}
M(t_0) = - \frac{6}{5}\hd^{\frac{3}{2}}\beta - \gamma\frac{3}{5}\frac{\pi\om^2}{\hd}\frac{\cos{\om t_0}}{\sinh{\frac{\pi\om}{\sqrt{\hd}}}}\left(\hd^2-\om^4\right).\end{equation*}
The energy balance between the linear damping and the excitation provided by mass variation which yields a zero of the
Melnikov Function occurs when
\begin{equation}
\gamma=\frac{2\hd^{\frac{5}{2}}\beta\sinh{\frac{\om\pi}{\sqrt{\hd}}}}{\om^2\pi\left(\om^4-\hd^2\right)\cos{\om t_0}}.\label{eq:gamma_F_M=0}
\end{equation}
For the parameters in Figure \ref{fig:Fgamma} this evaluates to $\gamma \approx .05681$, so for $\gamma$ above this
value the mass variation alone is enough to cause basin erosion. It is worth noting that this also occurs away from regions
of parametric instability. For $\om=\pi$ there is no parametric resonance near $\hd=1$, yet when $F=0$,
$\gamma_M \approx 6.467\beta$ still provides a zero of $M(0)$.

\subsubsection*{The General Effect of Time Varying Mass}

\begin{figure}[b]
\begin{center}
\includegraphics[width = 0.5\linewidth]{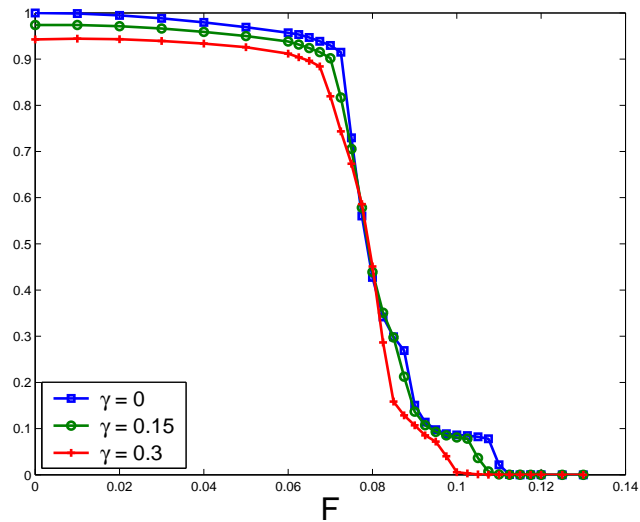}
\caption{Integrity Diagram for $\hd=\om^2=1$. The area of the safe
basin (i.e. the union of the basins of attraction of the systems
bounded attractors) normalised with respect to $F=\gamma=0$ is
plotted against the forcing parameter. }\label{fig:omdel1}
\end{center}
\end{figure}

\begin{figure}
\noindent
\begin{minipage}[b]{.43\linewidth}
\centering\includegraphics[width=\linewidth]{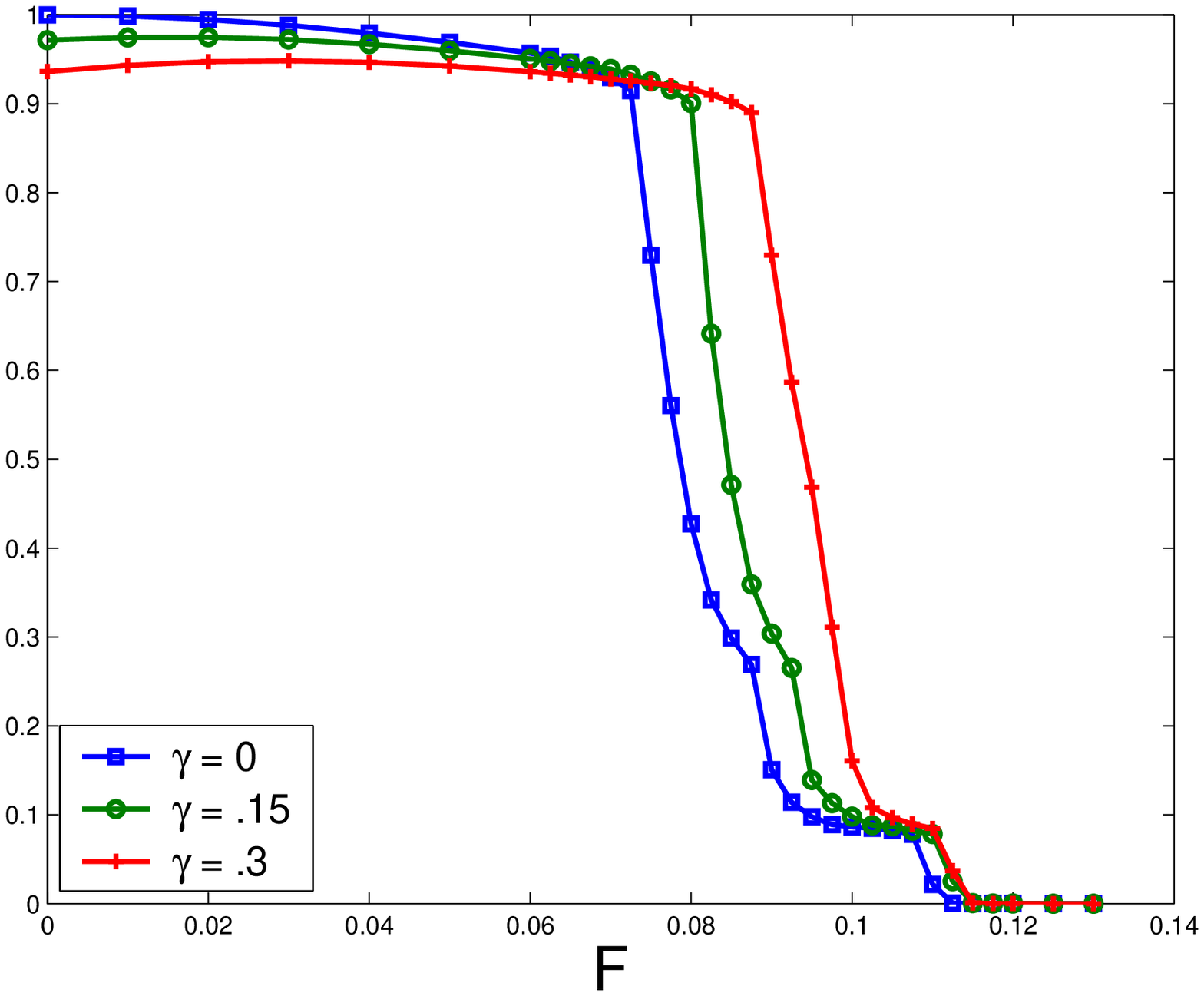}
\caption{$\hd = 1, \om = 0.85$}\label{fig:omdel2}
\end{minipage}
\begin{minipage}[b]{.43\linewidth}
\begin{center}
\includegraphics[width=\linewidth]{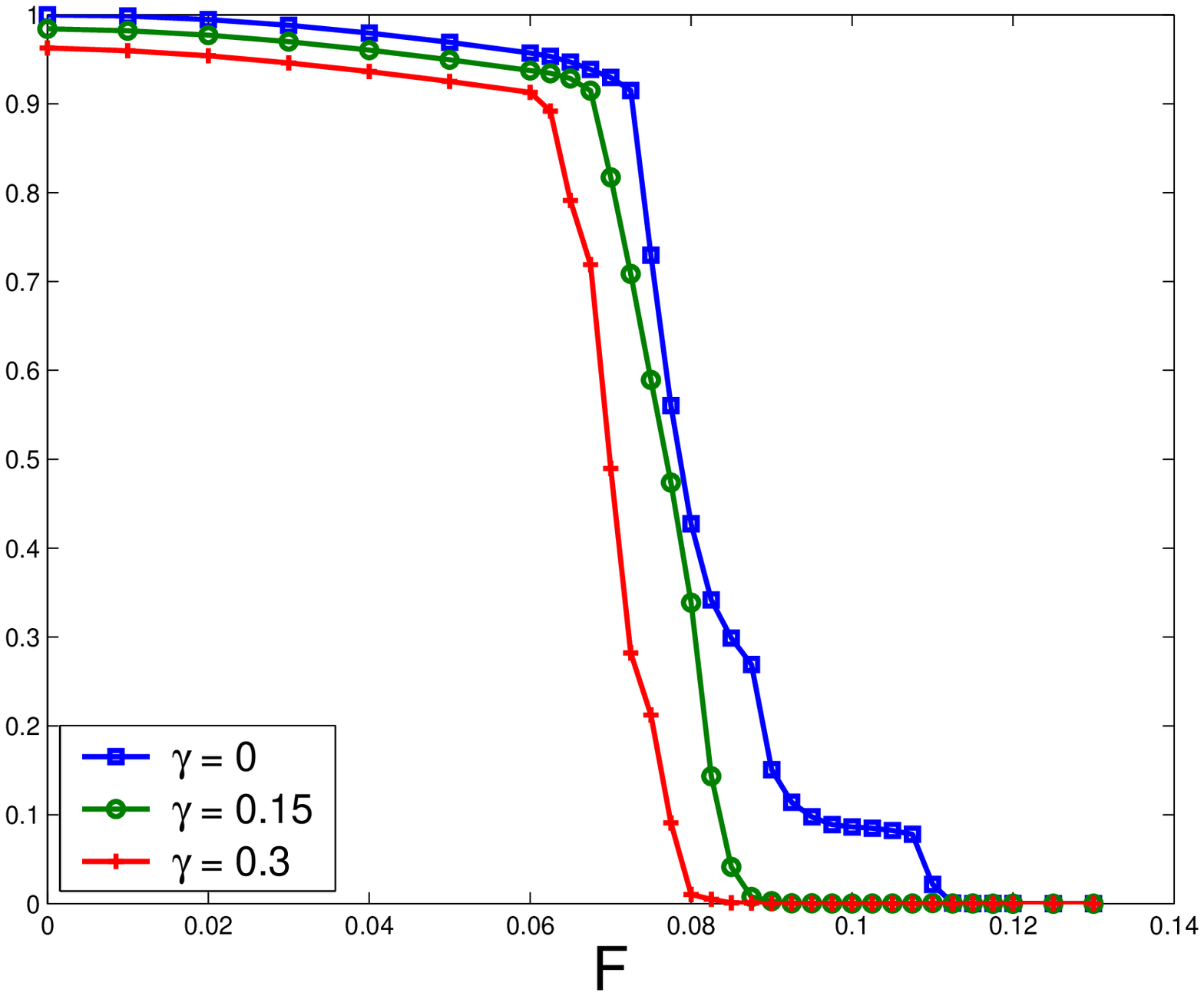}
\caption{$\hd=1, \om=1.15$}\label{fig:omdel3}
\end{center}
\end{minipage}
\end{figure}

The form of the Melnikov Function for periodic time varying mass \eqref{eq:Melnikov} can be considered as the usual
Helmholtz--Thompson Melnikov function plus a correction which takes into account the mass variation. The mass varying term
is
$$ - \gamma\frac{3}{5}\frac{\pi\om^2}{\hd}\frac{\cos{\om t_0}}{\sinh{\frac{\pi\om}{\sqrt{\hd}}}}\left(\hd^2-\om^4\right)$$
and this changes sign at $\hd =\pm\om^2$. This means that when
$\hd = \om^2$ the time varying mass does not affect the Melnikov
estimate of the first homoclinic tangle. Figure \ref{fig:omdel1}
is the {engineering integrity diagram} for this case. Although
there is some minor quantitative effect on the size of the initial
safe basin due to $\gamma>0$ the position of the dover cliff
reduction in basin area is not effected. If we set $M(t_0)=0$ and
solve for $F$ as usual we get
\begin{equation}F =\frac{\hd^{\frac{5}{2}}\beta}{5\pi\of^2} \frac{\sinh{\frac{\pi\of}{\sqrt{\hd}}}}{\cos{\of t_0}}
+\gamma\frac{\pi\om^2}{10\pi\of^2}\frac{\sinh{\frac{\pi\of}{\sqrt{\hd}}}}{\sinh{\frac{\pi\om}{\sqrt{\hd}}}}\frac{\cos{\om
t_0}}{\cos{\of t_0}}\left(\hd^2-\om^4\right).\end{equation}
Assuming $\gamma >0$ this makes it clear that when $\hd^2>\om^4$
the effect of the time varying mass is to delay basin erosion.
When $\hd^2<\om^4$ the basin erosion occurs for smaller $F$. This
effect is obvious in the integrity diagrams in Figures
\ref{fig:omdel2} and \ref{fig:omdel3}.
\section{CONCLUSIONS}
When Time Varying Mass in an oscillator is taken to be periodic the effect is similar to parametric forcing. Mathieu type
resonance occurs although due to coexistence only half as many regimes of parametric resonance exist. Also as the parameter
$\hd$ grows, i.e. as the mean mass of the oscillator tends to zero, the regimes of instability become thinner and thinner.
Whilst the linear parametric resonance does not have a direct qualitative effect on the dynamics of basin erosion, the
quantitative effects are significant and depend critically on the relationship between the mean mass and the frequency of
the variation. In some cases it can be broadly considered a stabilising influence whilst in others the opposite is true.

In a realistic model of a damaged ship the assumption of (sinusoidal) periodicity for the variation of the mass
due to water ingress is unlikely to
be realised exactly. However it is often the case that flood water mass has some periodic component. When this
occurs the effects studied in this paper are likely to be an important element of a complete picture of the dynamics.
\section*{ACKNOWLEDGEMENTS}
I would like to thank Professor Dracos Vassalos for suggesting
periodic time varying mass as an interesting problem to study and
providing support for my research. I would also like to acknowledge
the help and support of Dr. Michael Grinfeld and a number of helpful
discussions with Phillipos Samalekos.
\pagebreak

\bibliographystyle{plain}
\end{document}